\documentclass[10pt]{amsart}
\usepackage{amssymb,amsmath,amscd,amsfonts,amsthm,mathrsfs}
\usepackage{verbatim} 
\input xy
\xyoption{all}
\usepackage[all]{xy}  
\usepackage{color}

\usepackage[latin1]{inputenc}

\newcommand{\N}{\mathbb{N}}

\newcommand{\R}{\mathbb{R}}

\newcommand{\Ac}{\mathcal{A}}

\newcommand{\Cc}{\mathcal{C}}
\newcommand{\Dc}{\mathcal{D}}

\def\cchi{\raisebox{.45 ex}{$\chi$}}

\newcommand{\Ar}{\mathscr{A}}

\newcommand{\Hr}{\mathscr{H}}

\newcommand{\Mr}{\mathscr{M}}
\newcommand{\Nr}{\mathscr{N}}

\def\build#1_#2^#3{\mathrel{\mathop{\kern 0pt#1}\limits_{#2}^{#3}}}

\newtheorem{theorem}{Theorem}[section]
\newtheorem{proposition}{Proposition}[section]
\newtheorem{lemma}{Lemma}[section]
\newtheorem{remark}{Remark}[section]
\newtheorem{example}{Example}[section]
\newtheorem{corollary}{Corollary}[section]

\numberwithin{equation}{section}

\def \rmi{{\rm i}}

\title{Unboundedness of adjacency matrices of locally finite graphs}

\begin{document}
\author{Sylvain Gol\'enia}
\address{Mathematisches Institut der Universit\"at Erlangen-N\"urnberg,
Bismarckstr.\ 1 1/2 \\
91054 Erlangen, Germany}
\email{golenia@mi.uni-erlangen.de}
\date{Version of \today}
\begin{abstract}
Given a locally finite simple graph so that its degree is not bounded, 
every self-adjoint realization of the adjacency matrix is unbounded from above. 
In this note we give an optimal condition to ensure it is also
unbounded from below. We also consider the case of weighted graphs.
We discuss the question of self-adjoint extensions and prove an optimal criterium.
\end{abstract}

\maketitle
\section{Introduction}
The spectral theory of discrete Laplace operators and adjacency matrices acting on
graphs in useful for the study, among others, of some electrical networks, some gelling 
polymers and number theory, e.g.\ \cite{CDS, DS, DSV}. 
The study of random walk on graph is intimately linked with the study of the
heat equation associated to a discrete Laplace operator, see for instance 
\cite{Chu, MW}. In some recent papers \cite{Web, Woj, Woj2}, one works 
in the general context of locally finite graph and consider  a
non-negative discrete Laplacian,  see also \cite{KL, KL2} for
generalizations to Dirichlet forms. A key feature to obtain a Markov semi-group 
and to hope to apply these techniques is the boundedness from below (or from above) of a certain self-adjoint operator. In this note, we are interested in some self-adjoint realization of the adjacency matrix on locally finite graphs. We give some optimal conditions to ensure 
that the operator is unbounded from above and from below. 

We start with some definitions. Let $V$ be a countable set. Let
$E:=V\times V\rightarrow [0,\infty)$ and assume that
$E(x,y)=E(y,x)$, for all $x,y\in V$. We say that $G:=(E,V)$ is an
unoriented weighted graph with \emph{vertices} $V$ and
\emph{weights} $E$. In the setting of electrical networks, the weights correspond to the conductances. We say that $x,y\in V$ are \emph{neighbors} if
$E(x,y)\neq 0$ and denote it by $x\sim y$.  We say there is a
\emph{loop} in $x\in V$ if $E(x,x)\neq 0$.  
The set of \emph{neighbors} of 
$x\in E$ is denoted by $\Nr_G(x):=\{y\in E, E(x,y)\neq 0\}$. The
\emph{degree} of $x$ is by definition $d_G(x):=|\Nr_G(x)|$. 
A graph is \emph{connected}, if for all $x,y\in V$, there exists a
$x-y$ \emph{path}, i.e.,
there is a finite sequence  $(x_n)_{n=1,\ldots N+1 }$ in $V^{N+1}$ so that
$x_1=x$, $x_{N+1}=y$ and $x_n\sim x_{n+1}$, for all $n\in [1, N]$. Here $N$
denotes the \emph{length} of the path. When $x=y$, the path is called
a $x-$\emph{cycle}. A $x-$cycle of length $3$ is called a $x-$triangle.
When $E$ has its values in $\{0,1\}$ and when the
graph has no loop, the graph is called \emph{simple}. 
When $E$ is integer valued, it is a \emph{multigraph}.
We shall say that $E$ is \emph{bounded from below}, if $\inf\{E(x,y) |\, x,y\in
V$ and $E(x,y)\neq 0\}>0$. In the sequel, we suppose that $G$ is unoriented
weighted, has no loop and that $G$ is locally finite, i.e., $d_G(x)$
is finite for all $x\in V$. As a general rule, when no risk of
confusion arises, we drop the subscript $G$.  

We associate to $G$ the complex Hilbert space $\ell^2(V)$. We denote
by $\langle \cdot, \cdot\rangle$ and by $\|\cdot \|$ the scalar
product and the associated norm, respectively. 
By abuse of notation, we denote simply the space by $\ell^2(G)$. The
set of complex functions with compact support in $V$ is denoted by
$\Cc_c(G)$. There are several ways to define a Laplace
operator. One often considers the Laplacian defined by  
\begin{eqnarray}\label{e:defreg}
(\Delta_{G, \circ} f )(x) :=\sum_{x\sim y} E(x,y)\big(f(y)- f(x)\big), \mbox{
    with } f\in \Cc_c(G). 
\end{eqnarray}
See for instance \cite{Chu} for some other definitions. In this note,
we focus  on the analysis of the following off-diagonal Laplace
operator, the so-called \emph{adjacency   matrix}. We set: 
\begin{eqnarray}\label{e:def}
(\Ac_{G, \circ} f )(x) :=\sum_{x\sim y}E(x,y) f(y), \mbox{ with } f\in \Cc_c(G).
\end{eqnarray} 
Both of them are symmetric and thus closable. We denote the
closures by $\Delta_{G}$ and $\Ac_{G}$, their domains
by $\Dc(\Delta_G)$ and $\Dc(\Ac_G)$, and their adjoints by
$(\Delta_G)^*$ and  $(\Ac_G)^*$, respectively. 

In \cite{Woj2}, see also \cite{Jor, Web}, one shows that 
the operator $\Delta_G$ is essentially self-adjoint on
$\Cc_c(G)$, when the graph is simple. In particular, one has that
$\Delta_G=\Delta_G^*$. In contrast, even in the case of a simple graph
$G$, $\Ac_G$ may have many self-adjoint extensions, see \cite{MO, Mu}. 

We denote by $\eta_\pm(\Ac):=\dim\ker(\Ac^* \mp \rmi)\in
\N\cup\{+\infty\} $ the  \emph{deficiency indices} of the symmetric operator
$\Ac$. First, since the operator $\Ac$  commutes with the
complex conjugation, its deficency indices are equal, see
\cite{RS}[Theorem X.3]. We denote by $\eta(\Ac)$ the common value. 
Therefore $\Ac$ possesses some
self-adjoint extension. One has: $\Ac$ is essentially
self-adjoint on $\Cc_c(G)$ if and only if $\eta(\Ac)=0$.
Moreover, if $\eta(\Ac)$ is finite, the self-adjoint extensions can be
explicitly parametrized by the unitary group $U(n)$ in dimension
$n=\eta(\Ac)$.  In Remark \ref{r:inf} and in Proposition
\ref{p:OptiNelson}, we explain how to construct adjacency matrices
with deficiency indices $(+\infty, +\infty)$. Using the Nelson  
commutator Theorem, we prove in Section \ref{s:Nelson}:
\begin{proposition}\label{p:Nelson}
Suppose that  $\sup_x\max_{x\sim   y}|d(x)-d(y)|<\infty$ and
$\sup_x\max_{x\sim   y}|E(x)-E(y)|<\infty$, where $E(x):=\max_{y\sim
  x}E(x,y)$. Suppose also that $d$ or $E(\cdot, \cdot)$ is
bounded. Then $\Ac$ is 
essentially self-adjoint on $\Cc_c(G)$.  
\end{proposition} 
In Remark \ref{r:Be} and Proposition \ref{p:OptiNelson}, we prove the
optimality of these hypotheses. Given a finite sequence of graph $G_n$
such that $\Ac_{G_n}$ is essentially self-adjoint on $\Cc_c(G_n)$, one
can consider the direct sum $\Ac_{G}:=\oplus_{i=1, \ldots, n}
\Ac_{G_i}$ defined on $G:= \cup_{i=1, \ldots, n} G_n$ and infers that
$\Ac_{G}$ is essentially self-adjoint on $\Cc_c(G)$. Using the
Kato-Rellich lemma, it is of common knowledge that the result 
remains true if one perturbs the
structure of the graph $G$ on a finite set, as the perturbation is
of finite rank. In Lemma \ref{l:surgery} and under some conditions, we
explain how to extend this result for an infinite sum and to a 
perturbation with support on a non-finite set, see also Corollary
\ref{c:ess}. In analogy to similar constructions on manifolds, we call
this procedure \emph{surgery}.   

We turn to the main interest of this note, the unboundedness of the
self-adjoint realizations of the adjacency matrix. It is well-known 
that the operator $\Ac_G$ is bounded if $d_G$ and $E(\cdot, \cdot)$ are bounded. The reciprocal is true if $E$ is bounded from below, see Proposition \ref{p:estbd} for a refined statement. 

The first statement is easy and will be proved in Section \ref{s:gene}:  

\begin{proposition}\label{p:main}
Let $G=(E,V)$ be a locally finite graph. Let
$\hat \Ac$ be a self-adjoint realization of the $\Ac$. If the weight $E$ is unbounded, 
then  $\hat \Ac$ is unbounded from above and from below. 
\end{proposition}

We now deal with bounded weights $E$ and will restrict to the case $E$ bounded from below in the introduction. We refer to Section \ref{s:gene} for more general statements.
Suppose also that $d_G$ is unbounded. Let
$\kappa_d(G)$ be the filter generated by $\{x\in V, d_G(x)\geq n\}$,
with $n\in \N$. We introduce the \emph{lower local complexity} of a graph $G$ by: 
\begin{align}\label{e:lim}
C_{\rm loc}(G)&:=\liminf_{x \rightarrow \kappa_d(G)}
\frac{N_G(x)}{d^2_G(x)}, \mbox{ where }  
N_G(x):=\left|\{x-\mbox{triangles}\}\right|,\\
\nonumber
&:= \inf \bigcap\left\{ \overline{\left\{\frac{N_G(x)}{d^2_G(x)},\, x\in V \mbox{ and }d_G(x)\geq n\right\}}, n\in \N\right\}.
\end{align}  
Here $x\rightarrow \kappa_d(G)$ means converging to infinity along the
filter $\kappa_d(G)$. Recall that $G$ has no loop and beware that the
$x-$triangle given by $(x, y, z, x)$ is different from the one given
by $(x, z, y, x)$. In other words a $x-$triangle is oriented.

We introduce also the refined quantity, the \emph{sub-lower local complexity} of a graph $G$:
\begin{align}\label{e:lim2}
C_{\rm loc}^{\rm sub}(G)&:= \inf_{\{G'\subset G,\, \sup d_{G'}=\infty\}}  C_{\rm loc}(G'),
\end{align}  
where the inclusion of weighted graph is understood in the following sense:
\begin{eqnarray}\label{e:inclu}
G'=(E', V')\subset G \mbox{ if } V'\subset V \mbox{ and } E':=E|_{V'\times V'}.
\end{eqnarray}  
This means we can remove vertices but not edges. We conserve the induced weight.
Easily, one gets:
\begin{eqnarray}\label{e:incluorder}
0\leq C_{\rm loc}^{\rm sub}(G)\leq C_{\rm loc}(G)\leq 1.
\end{eqnarray}  
The sub-lower local complexity gives an optimal condition to ensure the 
unboundedness, from above and from below, of the self-adjoint realizations of the adjacency matrix. We give the main result:
\begin{theorem}\label{t:main}
Let $G=(E,V)$ be a locally finite graph such that $d_G$ is unbounded. Let
$\hat \Ac$ be a self-adjoint realization of the $\Ac$. Suppose that 
$E$ is bounded. Then, one has:
\begin{enumerate}
\item $\hat \Ac$ is unbounded from above. 
\item If  $C_{\rm loc}^{\rm sub}(G)=0$  and $E$ is bounded from below. Then $\hat \Ac$ is unbounded from below.
\item For all $\varepsilon>0$, there is a connected simple graph $G$
  such that $C_{\rm  loc}(G)\in (0, \varepsilon)$, $\Ac$ is essentially
  self-adjoint on $\Cc_c(G)$  and is bounded from below.
\end{enumerate}  
\end{theorem} 
By contrast, for any locally finite graph, $\Delta$ is a non-negative
operator, i.e.\ $\langle f, \Delta f\rangle\geq 0$, for all $f\in
\Dc(\Delta)$. 
The two first points come rather easily, see Proposition
\ref{p:nonbd} for a more general statement. 
The main difficulty is to prove the optimality given in the
last point, see Section \ref{s:opt}. 

\begin{example}
Consider $G$ a simple graph. If a graph $G$ has a subgraph, in the sense of \eqref{e:inclu}, being $\cup_{n\geq 0} S_{u_n}$ for some sequence $(u_n)_{n\in \N}$ that tends to infinity, 
then $C_{\rm loc}^{\rm sub}(G)=0$.
Here, $S_n=(E_n, V_n)$ denotes the \emph{star graph} of order $n$, i.e., $|V_n|=n$ and there is $x_\circ\in V_n$ so that $E(x, x_\circ)= 1$ for all $x\neq x_\circ$ and $E(x,y)=0$ for all $x\neq x_\circ$ and $y\neq x_\circ$. 
\end{example}
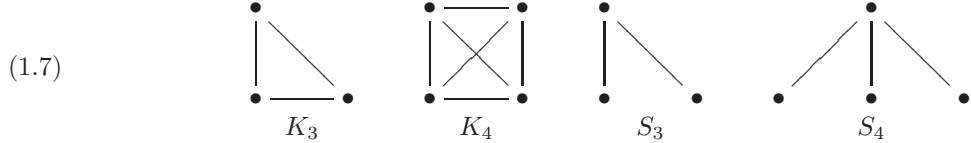
\begin{align}\label{e:SK}
\begin{array}{ccccccc}
\xymatrix{
{\bullet}\ar@{-}[d]\ar@{-}[dr]&
\\
{\bullet}\ar@{-}[r]& {\bullet}
} 
&&
\xymatrix{
{\bullet} \ar@{-}[r] \ar@{-}[d] \ar@{-}[dr]&{\bullet} \ar@{-}[d] \ar@{-}[ld]
\\
{\bullet}\ar@{-}[r]&{\bullet}}
&&
\xymatrix{
{\bullet}\ar@{-}[d]\ar@{-}[dr]&
\\
{\bullet}& {\bullet}
} 
&&
\xymatrix{
&{\bullet}\ar@{-}[d]\ar@{-}[dr]\ar@{-}[dl]&
\\
{\bullet}&{\bullet}& {\bullet}
} 
\\
K_3
&&K_4
&&S_3
&&S_4
\end{array}
\end{align}

We recall the definition of $K_n:=(E_n, V_n)$ the \emph{complete graph} of $n$ elements: $V_n$ is a set of $n$ elements and $E(a,b)=1$ for all $a, b\in V_n$, so that
$a\neq b$. One has $N_{K_n}(x)/d^2_{K_n}(x)= (n-1)(n-2)/n^2$, for all $x\in V_n$. 
Therefore, one can hope to increase the lower local complexity by having a lot of complete graphs as sub-graph in the sense of \eqref{e:inclu}.
More precisely, it is possible that $C_{\rm loc}(G)$ is positive, whereas $C_{\rm loc}^{\rm sub}(G)=0$. For instance, one has:

\begin{example}
For all $\alpha\in \N^*$, there is a simple graph $G$ such that 
$0=C_{\rm loc}^{\rm sub}(G)<C_{\rm loc}(G)= 1/(1+\alpha)^2$.
Now we construct the graph $S_{m+1} K_n=(E_{m,n},V_{m,n})$ as follows. Take $V_{m,n}:=\{x_n, x_{n,1},\ldots, x_{n, m+n}\}$. Set $E_{m,n}(x_n, x_{n,j})=1$, for all $j=1, \ldots, m+n$, $E_{m,n}(x_{n, j}, x_{n, k})=0$, for all  $j, k=1, \ldots, m$, and $E_{m,n}(x_{n, j}, x_{n, k})=1$, for all  $j, k=m+1, \ldots, m+n$, with $j\neq k$.
Set $G_\circ:=(E_\circ, V_\circ)$ as $\cup_{n\in \N^*} S_{\alpha n+1} K_n$. 
Finally, consider $G:=(E, V)$, with $V:=V_\circ$ and $E(x,y):=E_\circ(x,y)+ \sum_{n\in \N^*}\delta_{\{x_n\}}(x)\delta_{\{x_{n+1\}}}(y)$ for all $x,y\in V$, where $\delta$ is the 
Kronecker delta.
\end{example}
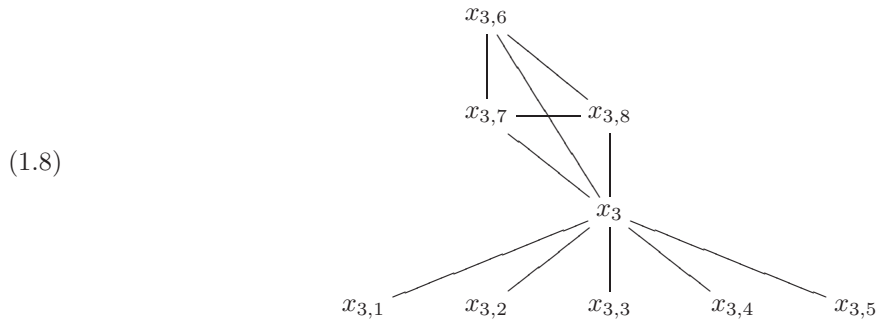
\begin{align}\label{e:SK3}
\begin{array}{ccc}
\xymatrix{
&{x_{3,6}}\ar@{-}[d]\ar@{-}[dr]&&&
\\
&{x_{3,7}}\ar@{-}[r]& {x_{3,8}}&&
\\
&&{x_3}\ar@{-}[dl]\ar@{-}[lld]\ar@{-}[d]
\ar@{-}[dr]\ar@{-}[drr]\ar@{-}[u]\ar@{-}[ul]\ar@{-}[uul]&&
\\
{x_{3,1}}& {x_{3,2}}&{x_{3,3}}&{x_{3,4}}& {x_{3,5}}
} 
\end{array}
\end{align}

We mention also that the (sub-)lower local complexity does not imply the essential self-adjointness of the adjacency matrix, see for instance Example \ref{ex:tree} and Remark \ref{r:Mu}. We point out that we know no example of a simple graph having the properties that the sub-lower local complexity is non-zero and that a self-adjoint realization of the adjacency matrix is unbounded from below.

In the Section \ref{s:Nelson} and in Section \ref{s:app}, we give some
criteria of essential self-adjointness for the adjacency matrix. 
In Section \ref{s:Nelsonopt}, we prove the optimality of the former
criterium. In Section \ref{s:gene}, we prove the Proposition and the
first part of the Theorem. In Section \ref{s:opt}, we prove the
optimality of the result by constructing a series of graphs and by
proceeding by surgery. At last in Section \ref{s:app}, we use the
surgery to give some examples with infinite deficiency indices. 

\noindent{\bf Notation:} We denote by $\N$ the set of non-negative
integers and by $\N^*$ the one of positive integers.
\\
\noindent{\bf Acknowledgments:} We would like to thank Thierry Jecko,
Andreas Knauf, and Hermann Schulz-Baldes for helpful discussions and
also grateful to  Daniel Lenz for valuable comments on the
manuscript. We warmly thank Bojan Mohar for having sent us his
reprint.

\section{Self-adjointness of the adjacency matrix}\label{s:adj}
One has the obvious inclusion $\Dc(\Ac)\subset \Dc(\Ac^*)$. 
The operator $\Ac$ is essentially self-adjoint on $\Cc_c(G)$ if 
$\Dc(\Ac)=\Dc(\Ac^*)$. In this case, there exists only one
self-adjoint operator $H$ so that $\Ac \subset H$, where the
inclusion is understood in the graph sense. In general, given a
self-adjoint extension $H$ of $\Ac$, one has: $\Ac \subset H=H^*
\subset \Ac^*$ and $\dim \big(\Dc(A^*)/\Dc(A)\big)=2\eta(A)$. 
The domain of $\Ac^*$ is given by:
\begin{eqnarray*}
\Dc(\Ac^*)=\Big\{f\in \ell^2(G), x\mapsto \sum_{x\sim y}E(x,y)f(y) \in
\ell^2(G)\Big\} \mbox{ and } \Ac^* f(x)= \sum_{y\sim x} E(x,y)f(y), \mbox{ for } f\in \Dc(\Ac^*).
\end{eqnarray*} 
It is well-known that the adjacency matrix of a locally finite graph
$G=(E,V)$ is usually not essentially self-adjoint. Using Jacobi
matrices it is easy to construct an example of weighted graphs with
$\max_x(d(x))\leq 2$, see Remark \ref{r:Be}. The first examples of
simple graphs are independently due to \cite{MO, Mu}. 

\subsection{A Nelson criterium}\label{s:Nelson} Under the hypothesis
of Theorem \ref{t:main}, one sees there is \emph{a   priori} no
canonical extension for $\Ac$, such as the Friedrich extension. 
Using the Nelson commutator theorem, we prove the criterium of
essential self-adjointness for $\Ac$ stated in the introduction.

\proof[Proof of Proposition \ref{p:Nelson}] Take $f\in \Cc_c(G)$. 
For $d$  bounded consider $\Mr(x):=E(x)$ and $\Mr(x):=d(x)$ when $E$
is bounded. Let $\Mr$ be the operator of multiplication by
$\Mr(\cdot)$. Then for some $c,C>0$, independent from $f$, we have: 
\begin{align*}
\|\Ac f\|^2&= \sum_x |\sum_{y\sim x} E(x,y)f(y)|^2\leq \sum_x
d(x)E^2(x)\sum_{y\sim   x} |f(y)|^2\leq \sum_x d(x) \max_{y\sim x}(d(y))E^2(x)
|f(x)|^2 
 \\
&\leq  \sum_x E^2(x)d(x)\big(c+d(x)\big) |f(x)|^2 \leq C\|\Mr f\|^2.
\end{align*}
Moreover:
\begin{align*}
|\langle f, [\Ac, \Mr] f\rangle|&= \left| \sum_x \overline{f(x)} 
\sum_{y\sim x} E(x,y)\big(\Mr(y)-\Mr(x)\big) f(y)\right|
\leq \sum_x \sum_{y\sim x} c|E^{1/2}(x)f(x)|\, |E^{1/2}(y)f(y)|.
\\
&\leq c \sum_x d(x)  |E^{1/2}(x)f(x)|^2 \leq
C \big\|\Mr^{1/2} f\big\|^2.
\end{align*}
Then using \cite{RS}[Theorem X.36], the result follows.\qed

\subsection{Optimality}\label{s:Nelsonopt} We now discuss the optimality
of the condition given in Proposition \ref{p:Nelson}. We start with a
remark.  

\begin{remark}\label{r:Be}
When $d$ is bounded, the condition on $E$ is optimal. Indeed,
set $\alpha>0$ and consider the Jacobi matrix acting on $\Cc_c(\N^*)$
with $0$ on the diagonal and $n^{1+\alpha}$, with $n\in \N^*$, on the
upper and lower diagonals. Then this adjacency matrix has deficiency
indices $(1,1)$, c.f., \cite[page 507]{Be} for instance. Here one has
$\max_{m=n\pm 1}|E(n)-E(m)|$ is equivalent to $n^\alpha$, when $n$ goes to infinity.
In this example, one can describe all the self-adjoint extensions by adding a condition at
infinity, using the Weyl theory, see \cite{SB} for recent results in
this direction. 
\end{remark} 

We now mimic the example of \cite{MO} in order to prove the optimality
when $E$ is bounded. 

\begin{proposition}\label{p:OptiNelson}
Let $\alpha>0$. There are $M\geq 1$ and a connected simple graph $G=(E,V)$, with
$V=\N$ and $n\sim n+1$ for all $n\in \N$, 
so that  $n^{-\alpha}|d(n)-d(n+1)|\leq M$, for all $n\in 
\N^*$ and so that $\eta(\Ac_G)=+\infty$. 
\end{proposition} 
\proof We construct a first graph $G^\circ=(E^\circ, V^\circ)$, with
$V^\circ:=\N$. It is a tree. Let $F:x\mapsto (x+1)^{\alpha+2}$.  
Given $n<m$, we say that $m\sim n$ if one has $m\in \big[F(n), F(n+1)\big)$.
Note that given $n\in \N^*$, there is a unique neighbor $m<n$ so that
$m\sim n$. We denote it by $n'$.  
Consider $f$ a solution of $\Ac_{G^\circ}^* f= \rmi f$. By
convention, set $f(0')=0$. We get: 
\begin{eqnarray}\label{e:indu}
f(n')+ \sum_{k\in [F(n), F(n+1))} f(k) = \rmi f(n), \mbox{ for }
  n\in \N.
\end{eqnarray} 
We construct $f$ inductively. We take $f(0)\neq 0$. Now, by using
\eqref{e:indu}, we choose  $f(k)$ constant, on the interval $\big[F(n),
F(n+1)\big)$. We denote by $c_n$ the common value and get:
\begin{eqnarray*}
c_n= \frac{1}{d(n)-1}\big(\rmi f(n)- f(n')\big),  \mbox{ for }
  n\in \N^*.
\end{eqnarray*} 
and $c_0= \rmi f(1)/d(0)$. Easily, $f$ is bounded. Moreover,
we have: 
\begin{eqnarray*}
\sum_{n=1}^M |f(n)|^2\leq 2 \sum_{n=1}^{F^{-1}(M)} \frac{|f(n)|^2+
  |f(n')|^2}{d(n)-1} \leq 4 \|f\|_\infty^2 \sum_{n=1}^\infty
\frac{1}{d(n)-1}.  
\end{eqnarray*} 
Note that $d(n)-1$ is equivalent to $(\alpha+2) (n+1)^{\alpha+1}$,
when $n$ goes to infinity. The series
converges and $f\in \ell^2(\N)$. Now remark that
as long as \eqref{e:indu} is fulfilled, it is enough to prescribe that $f(k)$
is constant on $\big[F(n), F(n+1)\big)$ for $k$ big enough. We infer the
deficiency indices of $\Ac_{G^\circ}$ are   infinite.   
Note also that $\sup_{n\in \N^*} n^{-\alpha}|d(n)-d(n+1)|$ is finite, 
since $n^{-\alpha}|d(n)-d(n+1)|$ tends to $(\alpha+2)(\alpha+1)$. 

It remains to connect $n$ with $n+1$. We proceed in the spirit of
Lemma \ref{l:surgery}. We define $G=(E,V)$ with $V=\N$. We say that
$m\sim n$ if $E^\circ(m,n)\neq 0$, for $m,n\in \N$ or if $|m-n|=1$, for
$m,n \in \N^*$. Now remark that $\|(\Ac_{G^\circ} -\Ac_G)f\|\leq
2\|f\|$, for all $f\in \Cc_c(G)$.  We recall that for a general
symmetric operator $H$, we have the topological direct sum $\Dc(H^*)=
\Dc(H)\oplus \ker(H^*+\rmi)\oplus \ker(H^*-\rmi)$.   To conclude, note
that $\Dc(\Ac_G)=\Dc(\Ac_{G^\circ})$ and $\Dc(\Ac_G^*)=\Dc(\Ac_{G^\circ}^*)$. \qed

\section{Unboundedness properties}
\subsection{Criterium of unboundedness}\label{s:gene}
In this section we give some elementary properties of the operator
$\Ac_G$ defined on a locally finite graph $G=(E,V)$.
We recall that if $E$ bounded from below, by definition, 
there is a $E_{\min}>0$, so that $E$ is with values in $\{0\}\cup[E_{\min}, \infty)$. 

To our knowledge, it is an open problem to characterize exactly the boundedness of the adjacency matrix of a graph with the help of the degree and the weights. For the Laplacian, one can show that it is bounded if and only if $\sup_x \big(\sum_{y\sim x} E(x,y)\big)$ is finite, e.g., \cite{KL2}. One has:
 
\begin{proposition}\label{p:estbd}
Let $G=(E,V)$ be a locally finite graph. Let $\hat \Ac_G$ be a self-adjoint realization of $\Ac_G$. Thus, one obtains:
\begin{eqnarray*}
\sup_{x\in V} \sum_{y\sim x} E^2(x,y)\,\leq \, \sup\sigma(\hat \Ac_G^2)\,\leq\, \sup_{x\in V} \,\sum_{y\sim x} d(y)
  E^2(x,y),
\end{eqnarray*}
where the value $+\infty$ is allowed. In particular, 
assuming that $E$ is bounded from below, then $\hat\Ac_G$ is bounded is and only if $d$ and $E$ are bounded. In this case $\hat\Ac_G=\Ac_G$.
\end{proposition} 
\proof Take $f\in \Cc_c(G)$ and consider $x\in V$. For the second inequality, one has:
\begin{align*}
\|\hat\Ac_G f\|^2&= \sum_{x}|\sum_{y\sim x} E(x,y) f(y)|^2\leq \sum_{x}\sum_{y\sim x} d(x)E^2(x,y) |f(y)|^2
= \sum_{x} \left(\sum_{y\sim x} d(y)E^2(x,y)\right) |f(x)|^2.
\end{align*}
We consider now the first inequality and ask $f$ to be with non negative values. We infer:
\begin{align*}
\|\hat\Ac_G f\|^2&= \sum_{x}|\sum_{y\sim x} E(x,y) f(y)|^2\geq \sum_{x}\sum_{y\sim x} E^2(x,y) |f(y)|^2 = \sum_{x} \left(\sum_{y\sim x} E^2(x,y)\right) |f(x)|^2.
\end{align*}
We conclude  $\sup\big(\sigma(\hat \Ac_G^2)\big)\geq 
\sum_{y\sim x_0} E^2(x_0,y)$, by taking $f$ with support in some $x_0\in V$.
\qed

Unlike the Laplacian, the adjacency matrix is not non-negative. Thus, in order to analyze accurately the unboundedness of the later one should consider the unboundedness from above and from below. We will make an extensive use of the fact that, given a self-adjoint operator $H$ acting in a Hilbert space $\Hr$, one has
$\inf(\sigma(H))= \inf\{\langle f,H f \rangle, f\in \Dc(H)$ and $ \|f\|=1\}$.
We will also use it for $-H$. As we deal with subgraphs, we will add a subscript to 
the neighbors relation $\sim$ to emphasize the use of the subgraph structure.

\begin{proposition}\label{p:nonbd}
Let $G=(E,V)$ be a graph and $\hat \Ac_G$ be a self-adjoint
extension of $\Ac_G$. Then, 
\begin{enumerate}
\item  If $E$ is not bounded then, the
spectrum of $\hat \Ac$ is neither bounded from above nor from below. 
\item In the sense of inclusion of graphs \eqref{e:inclu}, one has:
\begin{eqnarray*}
\sup \sigma(\hat \Ac_G)\geq \sup_{G'\subset G}\sup_{x\in V(G')}\left( \frac{1}{\sqrt{d_{G'}(x)}}\sum_{y \build{\sim}_{G'}^{}x}E(x,y)
+ \frac{1}{2d_{G'}(x)}\sum_{y\build{\sim}_{G'}^{}x}\,\sum_{z\build{\sim}_{G'}^{}y, z\build{\sim}_{G'}^{}x} E(y,z)\right).
\end{eqnarray*}
In particular, if $d$ is not bounded and $E$ bounded from below, 
then the spectrum of $\hat \Ac_G$ is not bounded from above, 
\item Suppose there is $C>0$ so that $\inf \sigma(\hat \Ac_G)\geq -C$. Then, for all $G'\subset G$, in the sense of \eqref{e:inclu}, 
\begin{eqnarray}\label{e:nonbd31}
\frac{1}{C}\,\big(\sum_{y \build{\sim}_{G'}^{}x}E(x,y)\big)^2\leq\sum_{y\build{\sim}_{G'}^{}x}\,\sum_{z\build{\sim}_{G'}^{}y, z\build{\sim}_{G'}^{}x} E(y,z) + C d_{G'}(x),
\end{eqnarray}
for $x\in G'$. In particular, when $E$ is with values in $\{0\}\cup [E_{\min}, E_{\max}]$, with 
\,$0<E_{\min}\leq E_{\max}<\infty$. Recalling \eqref{e:lim} and \eqref{e:lim2}, one obtains:
\begin{eqnarray}\label{e:nonbd32}
\frac{1}{C}\frac{E_{\min}^2}{E_{\max}}\, \leq\,  C_{{\rm loc}}^{{\rm sub}}(G)\,\leq\, C_{{\rm loc}} (G).
\end{eqnarray}
\end{enumerate} 
\end{proposition} 
\proof Let $G'$ be a subgraph of $G$, in the sense of \eqref{e:inclu}. Fix $x\in V(G')$
and consider a real-valued function $f$ with support in $\{x\}\cup
\Nr_{G'}(x)$. We have 
\begin{align}\nonumber
\langle f, \hat \Ac_G f\rangle&= f(x) (\Ac_{G'} f)(x)+ \sum_{y\build{\sim}_{G'}^{} x}
f(y)(\Ac_{G'} f) (y) 
\\\label{e:nonbd}
&= 2 f(x) (\Ac_{G'} f)(x) + \sum_{y\build{\sim}_{G'}^{} x} f(y) \sum_{z\build{\sim}_{G'}^{}
  y, z \build{\sim}_{G'}^{} x} E(y,z) f(z).  
\end{align} 
We first consider the case. There is a sequence $(x_n, y_n)_n$ of elements of $V^2$, such that  $E(x_n, y_n)\rightarrow +\infty$, when $n$ goes to infinity. Take
$G'=G$ and $f=f_n$  with support in $\{x_n, y_n\}$ in \eqref{e:nonbd}. We get $\langle f_n,
\hat \Ac f_n\rangle= 2 E(x_n, y_n) f(x_n)f(y_n)$. Then, choose $f(y_n)=1$ and  
$f(x_n)=\pm 1$ and let $n$ tend to infinity.

For the second case, take $f(x)=1$ and $f(y)=d_{G'}(x)^{-1/2}$ for $y$ neighbor of $x$ in $G'$. Noting that $\|f\|^2=2$, \eqref{e:nonbd} establishes the result. 

Focus finally on the third point. Take $f(x)=1$ and $f(y)=b$ for $y$ neighbor of $x$ in $G'$. Note that $\|f\|^2=1+d_{G'}(x)b^2$. Now, since $\langle f, \hat \Ac f\rangle\geq -C \|f\|^2$, 
\eqref{e:nonbd} entails:
\begin{eqnarray*}
2 b \sum_{y\build{\sim}_{G'}^{} x} E(x,y)+ b^2\sum_{y\build{\sim}_{G'}^{} x} \sum_{z\build{\sim}_{G'}^{}  y, z \build{\sim}_{G'}^{} x} E(y,z)+C(1+d_{G'}(x)b^2)\geq 0,
\end{eqnarray*}
for all $b\in \R$. Thus, the discriminant of this polynomial in $b$ is non-positive. This gives directly \eqref{e:nonbd31}. In turn, this infers:
\begin{eqnarray*}
\frac{1}{C}\frac{E_{\rm min}^2}{E_{\rm max}}\leq\frac{N_{G'}(x)}{d^2_{G'}(x)}+ \frac{1}{E_{\rm max}}\frac{C}{d_{G'}(x)}
\end{eqnarray*}
The statement \eqref{e:nonbd32} follows right away by taking the limit inferior with respect to the filter $\kappa_d(G')$.\qed

We now construct a tree where the hypotheses of the previous
propositions are fulfilled.  
\begin{example}\label{ex:tree} Take $M\in \N^*$. Let $\Ar_n$ be a set of $n$
  elements.  A word of $K\in \N^*$ letters build out of the alphabet
  $\{\Ar_n\}_{n\in \N}$ with increment $M$ is an element of
  $\Ar_M\times\ldots \times \Ar_{KM}$. The word of  $0$ letter is the
  empty set. Let $V$ be this set of words. We say that $K(x)$ is the
  length of a word $x\in V$.  For $x,y\in V$, we say that $x\sim y$ if
  $|K(x)-K(y)|=1$ and if they are composed of the same $\max\big(K(x),
  K(y)\big)-1$ first letters. 
  Then the adjacency matrix  $\Ac$
  defined on $G=(E,V)$ is essentially self-adjoint on $\Cc_c(G)$ by
  Proposition \ref{p:Nelson} and unbounded from below and from above
  by Proposition \ref{p:nonbd}. 
\end{example}
\begin{align*}
\xymatrix{
&&&\emptyset\ar@{-}[dl]\ar@{-}[dr]\ar@{-}[drrrr]&&&&
\\
&&a_0\ar@{-}[d]\ar@{-}[dl]\ar@{-}[dll]&& a_1\ar@{-}[ld]\ar@{-}[d]\ar@{-}[dr]\ar@{-}[drr]&&&a_3\ar@{-}[d]
\\
a_0a_0\ar@{.}[d]&a_0a_1\ar@{.}[d]&a_0a_5\ar@{.}[d]&a_0a_1\ar@{.}[d]&a_0a_3\ar@{.}[d]&a_0a_4\ar@{.}[d]&a_0a_5\ar@{.}[d]&a_3a_5\ar@{.}[d]
\\
&&&&&&&
} 
\\
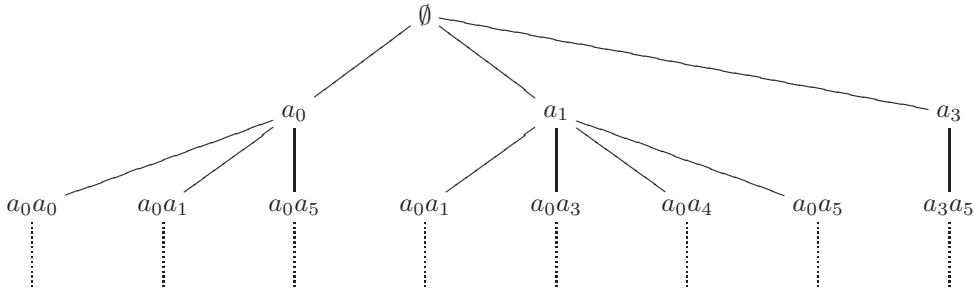
\mbox{Example of tree with increment } M=3. \mbox{\quad\quad\quad\quad\quad\quad\quad\quad\quad\quad}
\end{align*}

\begin{remark}\label{r:Mu}
In \cite{Mu}, one constructs a tree where the number of letters
increases exponentially and proves the adjacency matrix is not
essentially self-adjoint. In this context, Proposition \ref{p:nonbd}
yields that every self-adjoint extension is unbounded from below and
from above.
\end{remark} 

\subsection{Optimality}\label{s:opt}
We now show the optimality of the condition on the lower local complexity. 
For each $\varepsilon >0$, we find a graph $G=(E,V)$ and a sequence
$(x_n)_n$ of elements of $V$ such that $N(x_n)/d^2(x_n)$ tends to a
limit included in $(0, \varepsilon)$ and such that the adjacency
matrix associated to $G$ is bounded from below and essentially
self-adjoint on $\Cc_c(G)$. 

\begin{lemma}\label{l:Kn}
For each $k,n\in \N^*$, there is a finite graph $K_{k,n}$ and a point
$x_{k,n}\in K_{k,n}$ so that:
\begin{enumerate}
\item We have $\lim_{n\rightarrow \infty} N(x_{k,n})/d^2(x_{k,n})= 1/(2k^2)$.
\item The adjacency matrix $\Ac_{K_{k,n}}$ is bounded from below by $-4k$, in the form sense.
\end{enumerate} 
\end{lemma} 
\proof Consider first the graph given by the disjoint union
$K_{k,n}^\circ:=\{x_n\}\cup (K_{n})^{k}$, where $x_n$ is a point and 
$K_n:=(E_n, V_n)$ the complete graph of $n$ elements, i.e., $V_n$ is
a set of $n$ elements and $E(a,b)=1$ for all $a, b\in V_n$, so that
$a\neq b$, see \eqref{e:SK}. Then connect $x_n$ with each vertices of $(K_{n})^{k}$ to
obtain $K_{k,n}$. Note that $K_{1, n-1}=K_n$ and that the first point
is fulfilled. 

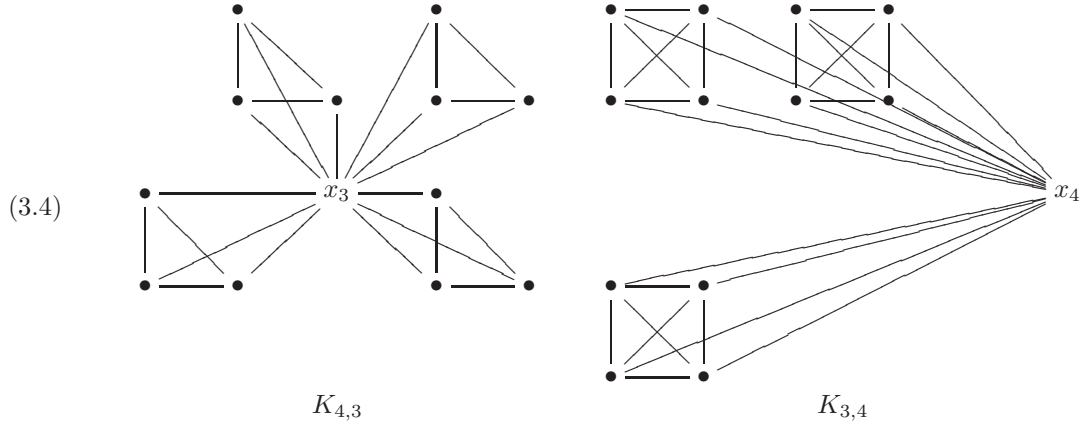
\begin{align}\label{e:SK2}
\begin{array}{ccc}
\xymatrix{
&{\bullet}\ar@{-}[d]\ar@{-}[dr]&&{\bullet}\ar@{-}[d]\ar@{-}[dr]&
\\
&{\bullet}\ar@{-}[r]& {\bullet}&{\bullet}\ar@{-}[r]& {\bullet}
\\
{\bullet}\ar@{-}[d]\ar@{-}[dr]&&{x_3}\ar@{-}[ru]\ar@{-}[rru]\ar@{-}[uur]\ar@{-}[dl]\ar@{-}[lld]\ar@{-}[ll]
\ar@{-}[r]\ar@{-}[dr]\ar@{-}[drr]\ar@{-}[u]\ar@{-}[ul]\ar@{-}[uul]&{\bullet}\ar@{-}[d]\ar@{-}[dr]&
\\
{\bullet}\ar@{-}[r]& {\bullet}&&{\bullet}\ar@{-}[r]& {\bullet}
} 
&&
\xymatrix{
{\bullet} \ar@{-}[r] \ar@{-}[d] \ar@{-}[dr]&{\bullet} \ar@{-}[d] \ar@{-}[ld]&{\bullet} \ar@{-}[r] \ar@{-}[d] \ar@{-}[dr]&{\bullet} \ar@{-}[d] \ar@{-}[ld]
\\
{\bullet}\ar@{-}[r]&{\bullet}&{\bullet}\ar@{-}[r]&{\bullet}
\\
&&&&&x_4 \ar@{-}[lluu] \ar@{-}[llluu] \ar@{-}[lllluu]\ar@{-}[llllluu]
\ar@{-}[llu] \ar@{-}[lllu] \ar@{-}[llllu]\ar@{-}[lllllu]
\ar@{-}[lllld]\ar@{-}[llllld]\ar@{-}[lllldd]\ar@{-}[llllldd]
\\
{\bullet} \ar@{-}[r] \ar@{-}[d] \ar@{-}[dr]&{\bullet} \ar@{-}[d] \ar@{-}[ld]&
&
\\
{\bullet}\ar@{-}[r]&{\bullet}&&
}
\\
K_{4,3}&&K_{3,4}
\end{array}
\end{align}
In a canonical basis, the adjacency  matrix of $K_{k,n}$
is represented by the $(n^k+1)\times(n^k+1)$ matrix:  
\begin{eqnarray*}
M(K_{k,n})=
\left(\begin{array}{ccccc}
0& 1 & 1& \cdots&1
\\
1& M(K_n) & 0&\cdots & 0
\\
1& 0 & M(K_n) & \cdots &0 
\\
\vdots &0&0& \ddots&0
\\
1& 0 & 0 & \cdots &M(K_n) 
\end{array}\right), \mbox{ where } M(K_n)=
\left(\begin{array}{ccccc}
0& 1 & 1& \cdots&1
\\
1& 0& 1&\cdots & 1
\\
1& 1 & 0 & \cdots &1 
\\
\vdots &1&1& \ddots& 1
\\
1& 1 & 1 & \cdots &0 
\end{array}\right).
\end{eqnarray*} 
Easily, the characteristic polynomial of $K_n$ is $\cchi_{K_n}(\lambda)=(-
\lambda +n-1)(-\lambda -1)^{n-1}$. Then we deduce that 
\begin{eqnarray}\label{e:poly}
 \cchi_{K_{k,n}}(\lambda)= \big(\lambda^2 - (n-1)\lambda -
 nk\big)(-\lambda+n-1)^{k-1}(-\lambda-1)^{k(n-1)}, 
\end{eqnarray} 
by replacing $C_1$ by $C_1- (\sum_{i\geq 2} C_i)/(-\lambda+n-1)$ in
the determinant of $M(K_{k,n})-\lambda$, for instance. Here $C_i$ denotes
the $i$-th column. At last, the second point follows from an  elementary
computation. \qed

\begin{remark}
We stress that there are only two subgraphs of $K_{k,n}$, the sense of
\eqref{e:inclu},  which are star graphs. Namely, $S_1$ and $S_k\simeq K_{k,1}$.
\end{remark} 

We now rely on a surgery lemma. 

\begin{lemma}\label{l:surgery}
Let $M\geq 1$. Given a sequence of graphs $G_n=(E_n, V_n)$, for
$n\in \N$. Choose $x_n\in V_n$. Let $G^\circ:=(E^\circ,
V^\circ):=\cup_{n\in \N} G_n$ be the disjoint union of $\{G_n\}_n$. Set $G:=(E, V)$
with $V=V^\circ$ and with $E(x,y):= E^\circ(x,y)$, when there is
$n\in \N$ so that $x,y\in V_n$ and where $\sup_{n\in \N} \sum_{m\in\N}
E(x_n, x_m) \leq M$.  
\begin{enumerate}
\item  We have $\|(\Ac_G -\Ac_{G^\circ})f\|\leq M
  \sup_{n,m}E(x_n, x_m)\|f\|$, for all $f\in   \Cc_c(G)=\Cc_c(G^\circ)$. 
\item The deficiency indices of  $\Ac_G$ are equal to
  $\eta(\Ac_G)=\sum_{n\in   \N} \eta(\Ac_{G_n})$. 
\item In particular, if $G_n$ are all finite graphs then $\Ac_G$ is
  essentially self-adjoint on $\Cc_c(G)$. 
\end{enumerate} 
\end{lemma} 
\proof We start with the first point. Observe that each $x_m$ has at
most $M$ neighbors in $\{x_n\}_{n\in \N}$. Then, 
\begin{align*}
\|(\Ac_{G}-\Ac_{G^\circ})f \|^2&= \sum_{n\in N}
\left|\big((\Ac_{G}-\Ac_{G^\circ})f\big)(x_n)\right|^2
=  \sum_{n\in N} \big| \sum_{m\in \N\setminus\{n\}}
E(x_n, x_m) f(x_m)\big|^2
\\
&\leq M \sum_{n\in \N}  \sum_{m\in \N\setminus\{n\}} E^2(x_n, x_m)
|f(x_m)|^2  
\leq M^2 \sup_{n,m}E^2(x_n, x_m) \sum_{n\in \N} |f(x_n)|^2.
\end{align*}
We turn to the second point. As we have a disjoint union,
$\eta(\Ac_{G^\circ})= \sum_{n\in   \N} \eta(\Ac_{G_n})$. For a
general symmetric operator $H$, we have the topological direct sum
$\Dc(H^*)= \Dc(H)\oplus \ker(H^*+\rmi)\oplus \ker(H^*-\rmi)$.  
To conclude, note that $\Dc(\Ac_G)=\Dc(\Ac_{G^\circ})$
and  $\Dc(\Ac_G^*)=\Dc(\Ac_{G^\circ}^*)$ from the first point.\qed

Finally, we establish the main result. 

\proof[Proof of Theorem \ref{t:main}] The two first points are
proved in Proposition \ref{p:nonbd}. Consider the last one. Given
$\varepsilon>0$, we choose $k>\sqrt{1/2\varepsilon}$. Given $M=2$,
we apply the Lemma \ref{l:surgery} with $G_n:=K_{k,n}$, where the latter is
constructed in Lemma \ref{l:Kn} by taking $E(x_n, x_m)\in \{0,1\}$, in
order to make the graph connected. We obtain a graph $G$ such that
$\Ac_G$ is essentially self-adjoint on $\Cc_c(G)$ and so that
$\Ac_G\geq -4k-M$. \qed 

\subsection{Further applications of the surgery}\label{s:app} First, we
give another criterium of essential self-adjointness.  One can perturb
the graph on a compact set and keep the same property, by the
Kato-Rellich Lemma.

\begin{corollary}\label{c:ess}
Let $N\in\N$. Consider $N$ simple graphs $G_n$ of constant degree,
i.e.\ $E_{G_n}$ has its values in $\{0,1\}$ and $d_{G_n}$ is constant
on $G_n$. Then for any graph $G$ obtained by surgery, as explained in Lemma
\ref{l:surgery}, one has $\Ac_G$ is essentially self-adjoint on $\Cc_c(G)$.
\end{corollary} 
\proof By \cite{Woj2}, one has $\Delta_{G_n}$ is essentially
self-adjoint on $\Cc_c(G_n)$. Since the graphs $G_n$ are simple with
constant degree $d_{G_n}$, $\Ac_{G_n} = d_{G_n} - \Delta_{G_n}$ is
 essentially self-adjoint on $\Cc_c(G_n)$. Lemma \ref{l:surgery}
concludes. \qed

To finish, we can create some arbitrary and possibly infinite
deficiency indices and obtain 
the property to be unbounded from above and from below.

\begin{remark}\label{r:inf}
Consider a graph $G_0$, such that $\Ac_G$ is not essentially
self-adjoint and has deficiency indices $(k,k)$, see for instance
Remarks \ref{r:Be} and \ref{r:Mu}.   
Considering $N$ copies of $G$, where $N\in \N\cup\{+\infty\}$, and by
joining  each copy as in Lemma \ref{l:surgery}. Then the adjacency
matrix of the new graph $G_1$ has deficiency indices $(Nk,Nk)$. 
Consider now the star graph $S_n$. Using again  Lemma \ref{l:surgery} with
$G_1$ and $\cup_{n\geq 2} S_n$ to obtain a graph $G_2$. Theorem \ref{t:main} 
gives  that every self-adjoint realization of
$\Ac_{G_2}$ is unbounded from above and from  below and that the
deficiency indices are also $(Nk, Nk)$.  
\end{remark}


\begin{thebibliography}{xxxxxx} 
\bibitem[Ber]{Be} J.M.\ Berezanski\i: \emph{Expansions in eigenfunctions
  of selfadjoint operators},
Providence, R.I.: American Mathematical Society 1968. IX, 809
p.\ (1968). 
\bibitem[Chu]{Chu} F.R.K.\ Chung: \emph{Spectral graph theory}
Regional Conference Series in Mathematics.\ 92.\ Providence, RI:
American Mathematical Society (AMS).\ xi, 207 p.\  
\bibitem[CDS]{CDS} D.\ Cvetkovi\'c, M.\ Doob, and H.\ Sachs: \emph{Spectra of graphs. Theory and application}, Second edition. VEB Deutscher Verlag der Wissenschaften, Berlin, 1982.\ 368 pp.
\bibitem[DSV]{DSV} G.\ Davidoff, P.\ Sarnak, and A.\ Valette:
\emph{Elementary number theory, group theory, and Ramanujan graphs},
London Mathematical Society Student Texts, 55.\ Cambridge University Press, Cambridge, 2003. x+144 pp.
\bibitem[DS]{DS} P.G.\ Doyle and J.L.\ Snell:
\emph{Random walks and electric networks}, The Carus Mathematical Monographs, the Mathematical Association of America, 159 pp.\ (1984). 
\bibitem[Jor]{Jor} P.E.T.\ Jorgensen: \emph{Essential
  self-adjointness of the graph-Laplacian}  
J. Math. Phys. 49, No.\ 7, 073510, 33 p.\ (2008).
\bibitem[KL]{KL} M.\ Keller and D.\ Lenz: {\em Dirichlet forms and
  stochastic completeness of graphs and subgraphs}, preprint
  arXiv:0904.2985v1 [math.FA].
\bibitem[KL2]{KL2} M.\ Keller and D.\ Lenz: {\em Unbounded Laplacians on graphs:
Basic spectral properties and the heat equation},
  Math.\ Model.\ Nat.\ Phenom.\ Vol.\ 5, No.\ 2, 2009.
\bibitem[MO]{MO} B.\ Mohar and M.\ Omladic: \emph{The spectrum of
  infinite graphs with bounded vertex degrees}, Graphs, 
hypergraphs and applications, Teubner Texte 73 (Teubner, Leipzig,
1985), pp.\ 122--125.
\bibitem[MW]{MW} B.\ Mohar and W.\ Woess: {\em
A survey on spectra of infinite graphs}, 
J.\ Bull.\ Lond.\ Math.\ Soc.\ {\bf 21}, No.3, 209-234 (1989). 
\bibitem[M\"u]{Mu} V.\ M\"uller: \emph{On the spectrum of an infinite
  graph}, Linear Algebra Appl.\ 93 (1987) 187--189. 
\bibitem[RS]{RS}M.\ Reed and B.\ Simon: {\em Methods of Modern Mathematical   
Physics, Tome I--IV: Analysis of operators} Academic Press.
\bibitem[SB]{SB} H.\ Schulz-Baldes: {\em Geometry of Weyl theory for
  Jacobi matrices with matrix entries} preprint arXiv:0804.3746v1 [math-ph].
\bibitem[Web]{Web} A.\ Weber: {\em Analysis of the Laplacian and the
  heat flow on a locally finite graph}, arXiv:0801.0812v3 [math.SP]
\bibitem[Woj]{Woj} R.\ Wojciechowski:  {\em Heat kernel and
essential spectrum of infinite graphs}, Indiana Univ.\ Math.\ J.\ 58
  (2009), no.\ 3, 1419 -- 1442.
\bibitem[Woj2]{Woj2} R.\ Wojciechowski: {\em Stochastic
completeness of graphs}, Ph.D. Thesis, 2007, arXiv:0712.1570v2[math.SP].
\end{thebibliography}
\end{document}